\documentclass[a4paper,12pt]{article}
\usepackage[cp1251]{inputenc}
\usepackage[english,russian]{babel}
\usepackage{amsthm,amsmath,amssymb}
\usepackage{fullpage}

\title{\selectlanguage{russian} Эрнест Борисович Винберг\\
{\normalsize (1937--2020)}}
\author{Д.~В.~Алексеевский, М.~В.~Белолипецкий, С.~Г.~Гиндикин, \\ В.~Г.~Кац, Д.~И.~Панюшев, Д.~А.~Тимашев, \\ О.~В.~Шварцман, А.~Г.~Элашвили, О.~С.~Якимова}
\date{ 5 августа 2021 г.}

\usepackage{color}
\definecolor{darkblue}{rgb}{0,0,.5}
\definecolor{darkgreen}{rgb}{.2,0.5,.2}
\usepackage[colorlinks=true,linkcolor=magenta,urlcolor=darkgreen,citecolor=darkblue]{hyperref}

\newcommand{\ZZ}{\mathbb{Z}}
\newcommand{\RR}{\mathbb{R}}
\newcommand{\CC}{\mathbb{C}}
\newcommand{\LL}{\mathcal{L}}
\newcommand{\Aa}{\mathsf{A}}
\newcommand{\Gg}{\mathsf{G}}
\newcommand{\g}{\mathfrak{g}}
\newcommand{\n}{\mathfrak{n}}

\newcommand{\sgl}{\mathfrak{sl}}
\newcommand{\ad}{\operatorname{ad}}
\newcommand{\GL}{\operatorname{GL}}
\newcommand{\Hom}{\operatorname{Hom}}
\newcommand{\Hor}{\operatorname{Hor}}
\newcommand{\Img}{\operatorname{Im}}
\newcommand{\Or}{\operatorname{O}}
\newcommand{\SL}{\operatorname{SL}}
\newcommand{\SO}{\operatorname{SO}}
\newcommand{\Sp}{\operatorname{Sp}}
\newcommand{\Spin}{\operatorname{Spin}}
\newcommand{\tr}{\operatorname{tr}}
\newcommand{\Sym}{\mathcal{S}}
\newcommand{\Uni}{\mathcal{U}}
\newcommand{\ii}{\boldsymbol{i}}

\newtheorem*{theorem}{Теорема}
\newtheorem*{corollary}{Следствие}

\begin{document}

\maketitle

12 мая 2020 года ушёл из жизни выдающийся математик Эрнест Борисович Винберг.

Эрнест Борисович родился 26 июля 1937~г.\ в Москве. Его отец Борис Георгиевич Винберг работал инженером-электриком на заводе <<Динамо>>, а мать Вера Евгеньевна Похвальнова преподавала математику и физику, а затем работала инженером-расчётчиком. Во время войны семья Винберга была в эвакуации в Пензенской области, а в 1943~г.\ вернулась в Москву. С математикой Э.~Б.~Винберг познакомился и заинтересовался ею очень рано, и уже в средней школе твёрдо решил, что будет математиком. С шестого класса школы он посещал математические кружки в Московском университете на Моховой улице, успешно участвовал в математических олимпиадах и в 1954~г.\ поступил на механико-математический факультет МГУ.

Молодость Эрнеста Борисовича пришлась на <<лучезарное десятилетие>> мехмата МГУ (1960-е~гг.). Это было время расцвета отечественной математики и, в частности, математики в Московском университете. На мехмате работали учёные высочайшего уровня (А.~Н.~Колмогоров, И.~М.~Гельфанд, И.~Р.~Шафаревич и многие другие). Пройдя по этажам факультета в Главном здании МГУ на Воробьевых горах и побеседовав с сотрудниками, можно было получить консультацию по любому вопросу современной математики. В студенческой среде царил энтузиазм, жажда знаний и стремление к научному творчеству. Из студентов тех лет вышли многие выдающиеся математики, и курс, на котором учился Эрнест Борисович, был особенно сильным. 

Научным руководителем Эрнеста Борисовича в студенческие годы и в аспирантуре, куда Э.~Б.~Винберг поступил после окончания университета в 1959~г., был Евгений Борисович Дынкин. Для Эрнеста Борисовича этот выбор был исключительно важен, и они продолжали поддерживать близкие отношения до смерти Евгения Борисовича. Е.~Б.~Дынкин был не только выдающимся математиком, но и исключительно оригинальным педагогом. Несомненно, что это обстоятельство сильно повлияло на становление Эрнеста Борисовича как математика и педагога. Для него научные исследования и работа с учениками были одинаково важны.

Главным семинаром для Эрнеста Борисовича был семинар Е.~Б.~Дынкина по группам Ли, который начался в 1956-57 учебном году как чисто студенческий, но уже на следующий год там появились Ф.~А.~Березин, Ф.~И.~Карпелевич и несколько позже И.~И.~Пятецкий-Шапиро. Это  были состоявшиеся, блестящие молодые математики,  в студенческие годы участвовавшие в предыдущей версии  семинара Дынкина. Они охотно общались с младшими участниками. Из студентов постоянно ходили А.~А.~Кириллов, С.~Г.~Гиндикин и А.~Л.~Онищик, начинавший аспирантуру. К тому времени Дынкин прекратил активную работу в области групп Ли, но сохранил к ней интерес и был очень успешным руководителем семинара. Продолжением семинара Дынкина впоследствии стал известный семинар Винберга--Онищика по группам Ли и теории инвариантов.

В 1961~г.\ состоялся <<колмогоровский>> приём новых сотрудников на мехмат МГУ: по инициативе А.~Н.~Колмогорова и при поддержке ректора МГУ И.~Г.~Петровского было решено омолодить преподавательский состав, и новый набор был осуществлён в основном из числа однокурсников Э.~Б.~Винберга. Сам Эрнест Борисович получил приглашение от заведующего кафедрой алгебры А.~Г.~Куроша и, перейдя в заочную аспирантуру, был взят на кафедру ассистентом.

Но ещё за год до этого, будучи аспирантом 2-го года, Винберг проводил семинарские занятия для студентов первого курса по линейной алгебре. Он вёл их, по свидетельству очевидцев, мастерски и захватывающе интересно.
Аура, которую создавал Эрнест Борисович на занятиях, немало способствовала тому, что из студентов тех лет вышло
несколько известных математиков. Некоторые из студентов (В.~Г.~Кац, Б.~Н.~Кимельфельд, А.~Г.~Элашвили, а несколько ранее Д.~В.~Алексеевский и Б.~Ю.~Вейсфейлер) стали первыми учениками, а потом и соратниками Винберга в научной работе.

В 1962~г.\ Э.~Б.~Винберг защитил кандидатскую диссертацию, решив проблему об однородных выпуклых конусах и областях Зигеля, поставленную Пятецким-Шапиро (подробнее об этом см.\ ниже). В 1965~г.\ по предложению Куроша Эрнест Борисович прочёл свой первый (не считая спецкурсов) курс лекций по линейной алгебре и геометрии и в том же году получил должность доцента.

В начале 1970-х~гг.\ в академической карьере (но не в научном творчестве!) Эрнеста Борисовича 
наступил трудный период. Винберг был в числе подписавших 
<<письмо девяноста девяти>> в защиту математика А.~С.~Есенина-Вольпина, насильственно помещённого в психиатрическую больницу по политическим мотивам, после чего последовала реакция. 
Для Эрнеста Борисовича наступил длительный перерыв в чтении лекций и возникли препятствия к руководству аспирантами на мехмате, а докторскую диссертацию ему пришлось защищать дважды --- в 1971 и 1984 годах (обе работы были блестящими и относились к гиперболическим группам отражений --- см.\ ниже о работах Винберга в этой области).

Положение изменилось с началом перестройки. В 1990~г.\ Эрнест Борисович стал профессором кафедры алгебры мехмата МГУ и работал на кафедре, читая лекции и спецкурсы, ведя семинары и руководя научной работой своих учеников, до конца дней.

Научное творчество Э.~Б.~Винберга 
исключительно широко и разнообразно. 
Не менее замечательно, а, пожалуй, даже и уникально то, что интересы Эрнеста Борисовича оставались многообразными на протяжении всей его жизни. Изучив ещё одну область,
сделав её <<своей>> и получив в ней первоклассные результаты, он не забывал о ней, когда новые интересы отодвигали её на второй план. Наоборот, периодически возвращаясь к прежней тематике, он получал в ней всё новые и новые интересные результаты. Так, он постепенно освоил и включил в свой научный арсенал такие разделы алгебры и геометрии как:
\begin{quote}
риманова геометрия однородных пространств; теория однородных и квазиоднородных выпуклых конусов; группы и алгебры Ли над незамкнутыми полями; алгебраическая теория инвариантов, включающая в себя, в том числе, и разработанную им теорию орбит и инвариантов, ассоциированных с градуировками полупростых алгебр Ли (тэта-группы); гиперболические группы отражений и геометрия пространств Лобачевского; арифметические решётки;
сложность действий алгебраических групп и сферические многообразия;
слабо симметрические и коммутативные однородные пространства; коизотропные действия алгебраических групп; эквивариантная симплектическая геометрия; редуктивные алгебраические полугруппы;
трансцендентная теория инвариантов; неабелевы градуировки алгебр Ли (список неполон).
\end{quote}
Во всех этих разделах он получил выдающиеся результаты и был признанным экспертом; не
говоря уже просто о том, что он являлся создателем некоторых из них.

Наряду с отмеченными выше большими темами исследований Винберга, о которых речь пойдёт ниже, можно обнаружить и его отдельные неожиданные работы, выходящие далеко за эти рамки, например:
\begin{itemize}
  \item <<К теореме о бесконечномерности ассоциативной алгебры>> (№139 в списке публикаций Э.~Б.~Винберга на \verb|www.mathnet.ru|);
  \item <<О теореме Шёнфлиса--Бибербаха>> (№115 на \verb|www.mathnet.ru|);
  \item <<The two most algebraic K3 surfaces>> (№97 на \verb|www.mathnet.ru|);
  \item <<On some number-theoretic conjectures of V.~Arnold>> (№30 на \verb|www.mathnet.ru|).
\end{itemize}

До самых последних дней Эрнест Борисович успешно работал во многих областях и творчески общался со своими учениками и коллегами. Эрнест Борисович неоднократно помогал советами как математикам, так и людям иных специальностей, которые обращались к нему с математическими вопросами и всегда получали квалифицированный ответ по существу.

Говоря о влиянии печатных работ Э.~Б.~Винберга, хочется отметить, что все они представляют собой
тексты, которые \emph{можно понять} и из которых \emph{можно многому научиться}. В них присутствует неуловимое умение доступного и понятного изложения, 
основанное, несомненно, на том, что Эрнест Борисович \emph{понимает} (в высшем смысле) то, о чём пишет. И
потому интуитивно, как само собой разумеющееся, выбирает наиболее ясный, экономный и поучительный
способ изложения материала. В то же время, разнообразные примеры, интересные замечания и открытые проблемы, упоминаемые в тексте, дают понять заинтересованному читателю, что написанное представляет только вершину айсберга, и подсказывают темы для дальнейших раздумий и собственных исследований. Эти качества, в частности, проявились в блестящей серии энциклопедических обзоров в <<Итогах науки и техники>> ВИНИТИ, написанных Э.~Б.~Винбергом с соавторами в конце 1980-х~гг.\ (№№ 78, 80, 83--86 на \verb|www.mathnet.ru|).

Стоит также заметить, что (русский) язык и стиль статей Э.~Б.~Винберга был безукоризненным, и некоторые его
ученики с благодарностью вспоминают его советы по стилистике. В связи с этим
вспоминается следующий любопытный эпизод. Одна из самых известных и
цитируемых работ Винберга \cite{V86} возникла из-за опечатки в другой, не менее известной
работе~\cite{VKi78}. Эрнест Борисович пишет об этом в первом абзаце в~\cite{V86}. Сама фраза
с опечаткой выглядит так:
<<\emph{\dots\ если группа $B$ действует локально транзитивно на~$X$, то группа $B$ имеет лишь
конечное число орбит в~$X$}>>, и тут вместо второго $B$ имелось в виду~$G$.
Так вот, в устном докладе на семинаре Эрнест Борисович сказал, что
сама напечатанная фраза стилистически неаккуратна, и он никогда бы не смог так написать!
Если бы в этом месте имелось в виду то, что написано, то во второй раз вместо слов
<<\emph{группа~$B$}>> было бы сказано <<\emph{она}>>. 
Cама эта опечатка
осталась незамеченной авторами до той поры, пока несколько лет спустя Винберг не получил письмо
(тогда ещё не электронное) от одного математика, которому потребовалось именно это, более сильное и недоказанное в \cite{VKi78} утверждение, и он спросил у Эрнеста Борисовича о подробностях. Так родилась замечательная статья~\cite{V86}. 

Ведя активную исследовательскую работу и публикуя собственные результаты (около 150 научных работ в списке публикаций на \verb|www.mathnet.ru|), Э.~Б.~Винберг отдавал много сил редакционной деятельности, состоя в редколлегиях ряда ведущих математических журналов: <<Функциональный анализ и его приложения>> (с 2005~г.), <<Труды Московского математического общества>> (с 1999~г.), <<Journal of Lie Theory>> (с 1994~г.). Эрнест Борисович был бессменным главным редактором журнала <<Transformation Groups>> с момента его основания (1996~г.).

Эрнест Борисович Винберг был не только великим математиком, но и великим педагогом, Учителем. У него было много учеников. Замечательные годы ученичества под руководством Эрнеста Борисовича затем перерастали в научное сотрудничество, которое продолжалось до конца его дней.
За эти годы легко было убедиться в исключительности таланта Эрнеста Борисовича, как педагога и исследователя. Этот талант позволил ему воспитать плеяду первоклассных математиков и обнаруживать математические факты и пути их доказательства, которые никто кроме него обнаружить бы не смог.
Как педагога его отличала интуиция, терпение и исключительная ясность изложения проблем, которые он ставил своим ученикам.

Почти все ученики Э.~Б.~Винберга выросли на семинаре по группам Ли на мехмате МГУ, которым он руководил с 1961 года совместно с А.~Л.~Онищиком, а впоследствии и с другими своими учениками. Семинар Винберга--Онищика на протяжении более полувека был значительным событием московской математической жизни. Записки семинара легли в основу книги~\cite{VO}, которую отличает нестандартный и оригинальный подход в изложении теории групп и алгебр Ли, делающий её уникальным источником как для специалистов, так и для неофитов, приступающих к изучению этой области.

Помимо преподавательской работы на мехмате, Винберг в середине 1960-х годов недолгое время работал во 2-ой московской школе, когда И.~М.~Гельфанд и Е.~Б.~Дынкин с помощью своих младших коллег начинали налаживать там специализированное математическое образование (то, что создало феномен и легенду <<Второй школы>>).
В начале 1990-х~гг.\ Эрнест Борисович принимал деятельное участие в организации Независимого Московского университета 
и читал там базовые, но весьма нестандартные и глубокие курсы алгебры. Из записок этих лекций, изданных в НМУ в 1995~г., впоследствии вырос его <<Курс алгебры>> (№53 на \verb|www.mathnet.ru|), выдержавший много переизданий и ставший классическим университетским учебником. В последние 15 лет жизни Э.~Б.~Винберг был главным редактором журнала <<Математическое просвещение>>.

Эрнест Борисович был научным руководителем серии летних и зимних школ <<Алгебры Ли, алгебраические группы и теория инвариантов>>, проводившихся попеременно в Самарской области и в Москве во второй декаде XXI~в., первая из которых состоялась под Самарой в июне 2009~г., а последняя на данный момент, 8-я по счёту, --- в Москве в январе-феврале 2020~г. Для их организации много сделали ученик Винберга И.~В.~Аржанцев и А.~Н.~Панов, а сам Эрнест Борисович был <<душой>> этих школ, неизменным председателем программного комитета и постоянным лектором. Он придавал этим школам большое значение в деле распространения современных математических знаний, стараясь в своих лекционных курсах (2009, 2011, 2014, 2017~гг.) в доступной форме, с присущим ему мастерством, познакомить молодых слушателей с новейшими результатами в тех областях, которыми он занимался.

В 1997 году Э.~Б.~Винберг получил премию Александра фон Гумбольдта, на которую его выдвинул университет Билефельда. Это положило начало более чем 20-летнему периоду сотрудничества Винберга с математиками из Билефельда (Г.~Абельс, Й.~Меннике, Х.~Хеллинг) и регулярным летним визитам в Германию. Для Эрнеста Борисовича это составляло важную часть жизни, поскольку там, вдали от московской суеты, он мог полностью отдаться научному творчеству. Вместе с тем, будучи искренним патриотом, Эрнест Борисович никогда не рассматривал возможность эмиграции, не мысля свою жизнь вне родины.

Он был необыкновенно цельным человеком, очень интеллигентным и доброжелательным, но своими принципами не поступался никогда. С его уходом математическая жизнь уже не будет прежней.

Перейдём к обзору основных направлений научных исследований Эрнеста Борисовича Винберга и его вклада в математику.

\subsubsection*{Комплексные однородные области, выпуклые однородные конусы, однородные кэлеровы многообразия}

Эрнест Борисович в начале работы семинара Дынкина интересовался геометрическими работами Э.~Картана. Он сделал несколько обзорных докладов и получил первые результаты. В 1958 году его интересы решительно изменились после того как  И.~И.~Пятецкий-Шапиро, только вернувшийся в Москву из Калуги, начал регулярно рассказывать на семинаре об областях Зигеля --- многомерных обобщениях верхней полуплоскости, которые интересовали его в связи с многомерными автоморфными формами. К.~Зигель рассмотрел
пример области комплексных симметрических матриц с положительно определённой мнимой частью. Пятецкий-Шапиро называл  области в  $\CC^n$ вида $D(V)=\RR^n+\ii V$, где $V$ --- выпуклый конус без прямых, областями Зигеля 1-го рода (их называют также радиальными трубчатыми). Эти области не ограничены, но голоморфно эквивалентны ограниченным.

Общие многомерные автоморфные формы связаны с эрмитовыми симметрическими пространствами или, что эквивалентно, с симметрическими однородными комплексными
областями. Такие области называются классическими, если их группа автоморфизмов является классической. Кроме полуплоскостей Зигеля были известны ещё 3 серии классических областей, допускающих реализацию в виде областей Зигеля 1-го рода. Однако имеются классические области (например,  комплексный шар),  которые не допускают такой реализации, поскольку у областей Зигеля 1-го рода вещественная размерность остова границы совпадает с комплексной размерностью области.

И.~И.~Пятецкий-Шапиро предложил обобщение, которое уже работало для всех классических  областей.  Для выпуклого конуса $V\subset \RR^n$ пусть $F:\CC^m\times\CC^m\to\CC^n$ --- такая эрмитова полуторалинейная векторнозначная форма, что $F(w,w)$ лежит в замыкании конуса $V$ для любого $w\in \CC^m$ (\emph{$V$-положительная форма}).
Тогда область Зигеля 2-го рода $D(V,F)$ состоит из таких
пар $(z,w)\in\CC^n\times\CC^m$, что $$\Img(z)-F(w,w) \in V.$$

Если конус $V$ линейно однороден  с группой $G(V)$ и для каждого  $g\in G(V)$ существует такой элемент  $\lambda(g)\in\GL_m(\RR)$, что
$g(F(w,w))=F(\lambda(g)w,\lambda(g)w)$ (\emph{условия однородности}), то соответствующая
область Зигеля является аффинно однородной. Области 1-го рода отвечают случаю $m=0$.

Пятецкий-Шапиро нашёл реализацию всех классических областей как аффинно однородных областей Зигеля 2-го рода и применил их к изучению автоморфных форм  в этих областях. Сенсацией стало его решение проблемы Э.~Картана: одна из простейших аффинно однородных областей Зигеля в $\CC^4$ является не симметрической. Это открытие произошло практически на глазах участников семинара Дынкина.

Э.~Картан нашёл все однородные комплексные ограниченные области в размерностях 2 и 3, которые оказались симметрическими. Отдельно он описал симметрические области во всех размерностях.
В 1950-е годы был консенсус, что комплексные однородные области  всегда являются симметрическими. Это было доказано А.~Борелем и Ж.-Л.~Кошулем для простых групп автоморфизмов, а Дж.~Хано --- для
унимодулярных. Контрпример Пятецкого-Шапиро опроверг эту гипотезу.

В 1959~г.\ И.~И.~Пятецкий-Шапиро пригласил участвовать в проекте младших участников
семинара Э.~Б.~Винберга и С.~Г.~Гиндикина, сыграв тем самым исключительную роль в становлении Винберга как математика. Он предложил вначале Эрнесту Борисовичу понять природу выпуклых однородных конусов, которые появляются
в областях Зигеля для классических областей. Эти конусы являются выпуклыми, линейно однородными, не содержат прямых и являются самосопряжёнными для некоторых положительных билинейных форм. Винберг увидел, что классификация таких конусов эквивалентна классификации компактных йордановых алгебр, которая была известна. Несколько раньше эту связь рассматривал М.~Кёхер. Винберг ссылался на некоторые результаты Кёхера, но о его основной теореме он узнал только после сдачи заметки в печать.
Сегодня этот результат называется теоремой Кёхера--Винберга.

Эрнест Борисович показал, что существуют  4 серии конусов, связанных  с йордановыми алгебрами, которые в точности отвечают 4-м сериям классических областей. Есть ещё особый конус, связанный с октавными матрицами 3-го порядка. Он отвечает реализации особой 27-мерной комплексной симметрической области как области Зигеля 1-го рода.

Заметка Винберга о конусах содержала ещё один важнейший результат:
пример выпуклого линейно однородного, но не самосопряжённого, 5-мерного конуса (в меньших размерностях таких нет), что автоматически приводит к примеру несимметрической аффинно однородной области Зигеля 1-го рода в размeрности~5. Параллельно Винберг построил континуум неэквивалентных выпуклых линейно однородных конусов  в размерности~11.

Какое-то время казалось,что понятие областей Зигеля нуждается в расширении, чтобы охватить все комплексные однородные области. Однако в 1963~г.\ Э.~Б.~Винберг, С.~Г.~Гиндикин и И.~И.~Пятецкий-Шапиро опубликовали центральный результат теории:

\begin{theorem}[\cite{VGP63}]
Всякая комплексная однородная ограниченная область биголоморфно эквивалентна аффинно однородной области Зигеля 2-го рода.
\end{theorem}

Эта теорема была доказана чисто алгебраическими методами. На алгебре Ли группы автоморфизмов с отмеченной подалгеброй изотропии имеется дополнительная структура \emph{$j$-алгебры}, отвечающая почти комплексной структуре и кэлеровой форме. Условия интегрируемости комплексной структуры и точность кэлеровой формы дают аксиоматику $j$-алгебр. Основное утверждение об  однородных областях перефразируется в некоторое структурное утверждение о $j$-алгебрах, которое и доказывается в статье.


Эрнест Борисович построил замечательный вещественный аналог теории комплексных однородных ограниченных областей. Она включает теорию выпуклых однородных конусов, но шире. А именно, рассматриваются вещественные выпуклые аффинно однородные области без прямых. Имеется вещественный аналог комплексных областей Зигеля. Пусть  $V \subset \RR^n$, как и выше, выпуклый линейно однородный конус без прямых. По аналогии с $V$-положительными эрмитовыми векторнозначными формами рассматриваются $V$-положительные симметрические билинейные формы $F:\RR^m\times\RR^m\to\RR^n$, определяемые условием, что $F(u,u)$ лежит в замыкании конуса $V$ при любом $u\in\RR^m$. Тогда вещественная область Зигеля $D(V,F)$ состоит из таких пар $(x,u)\in\RR^n\times\RR^m$, что $$x-F(u,u)\in V.$$ Как и в комплексном случае, если выполняются условия однородности для некоторой транзитивной группы линейных преобразований конуса, то эта вещественная область Зигеля является аффинно однородной.

\begin{theorem}
Всякая аффинно однородная область без прямых в вещественном
пространстве аффинно эквивалентна аффинно однородной вещественной области Зигеля.
\end{theorem}

Возникает индуктивная конструкция вещественных аффинно однородных областей, включающая однородные конусы. Конус над $n$-мерной аффинно однородной вещественной областью Зигеля даёт выпуклый линейно однородный конус размерности $n+1$. И обратно, однородные вещественные области Зигеля можно  реализовать как гиперплоские сечения (параболического типа) выпуклых однородных конусов.

В частности, для однородных конусов  число индуктивных шагов можно интерпретировать как ранг (это согласуется с понятием ранга симметрического конуса). Соответственно, определяется картановская подгруппа, корни, корневые подпространства, кратности корней. Система корней однородного конуса всегда имеет тип~$\Aa$. Условие равенства кратностей имеется только в самосопряжённом случае, когда группой Вейля служит группа перестановок.

Центральное место в работе Винберга занимает введённый им  специальный класс алгебр, отвечающих выпуклым однородным конусам: \emph{$T$-алгебры}. Точки линейного пространства конуса реализуются как обобщённые матрицы, векторные элементы которых лежат в корневых подпространствах. При этом  аксиоматика умножений векторных матричных элементов такова, что  треугольные матрицы образуют группу. Эта группа действует на пространстве обобщённых матриц с открытой просто транзитивной орбитой (единичной матрицы), являющейся выпуклым однородным конусом.

Теория однородных выпуклых конусов и областей составила содержание кандидатской диссертации Э.~Б.~Винберга, см.~\cite{V63}. За эти результаты Винберг получил премию Московского математического общества (1963~г.,
совместно с С.~Г.~Гиндикиным)

Между комплексной и вещественной теориями  есть не только аналогия, но и более глубокая связь. Можно рассматривать  в $\CC^n$ трубчатые области $T(U)=\RR^n +\ii U$, где $U$ --- вещественная выпуклая область без прямых (не обязательно конус). Тогда область $T(U)$ будет аффинно однородной в комплексном смысле тогда и только тогда, когда $U$ будет аффинно однородной в вещественном смысле. Аффинно однородные трубчатые области могут быть реализованы
как трубчатые с аффинно однородной вещественной областью Зигеля.

При доказательстве основного результата о комплексных однородных ограниченных областях стало ясно, что решающую роль играет наличие на них структуры однородных кэлеровых многообразий.

На этом основании Э.~Б.~Винберг и С.~Г.~Гиндикин рассмотрели случай произвольных однородных кэлеровых многообразий. Кроме областей были известны 2 типа таких многообразий: локально плоские, которые получаются путём факторизации линейного эрмитова пространства по некоторой решётке, а также односвязные компактные, которые были полностью классифицированы. Возникла гипотеза, что произвольное однородное кэлерово многообразие допускает однородное расслоение с комплексной ограниченной однородной областью в качестве базы и прямыми произведениями локально плоских и односвязных компактных кэлеровых однородных многообразий в качестве слоев.

В 1967~г.\ Винберг и Гиндикин доказали гипотезу для центрального случая разрешимых групп автоморфизмов \cite{VG67}. Окончательно гипотеза была доказана в 1995~г.

Теория однородных выпуклых конусов не только занимает особое место в творчестве Э.~Б.~Винберга, но и находит удивительные применения далеко за её пределами. Оказалось, что эта теория допускает интерпретацию в рамках   информационной геометрии Ченцова--Амари, как геометрия одного из самых интересных классов статистических  многообразий, так называемых \emph{экспоненциальных семейств}. 
Согласно С.~Амари, экспоненциальное семейство статистических распределений <<не только является типичной статистической моделью, включающей  многие хорошо известные  семейства  вероятностных распределений, такие, как дискретное распределение, распределение Гаусса, мультинормальное распределение, гамма-распределение и т.\,д., но и ассоциируется с выпуклой функцией, известной как кумулятивная порождающая  функция или свободная энергия>>~\cite{Am}.

Информационная геометрия  изучает  многообразия  вероятностных распределений, снабжённые дивергенцией  (несимметричным расстоянием)  и  ассоциированной  с ней  римановой метрикой (\emph{метрикой Фишера--Раo}). С выпуклым конусом $V\subset\RR^n$ связана характеристическая функция Винберга--Кошуля $\varphi_V:V\to\RR_{>0}$ (ранее также возникавшая у Бохнера в связи с трубчатыми областями и преобразованием Лапласа):
$$\varphi_V(x)=\int_{V^*}e^{-\langle x,y\rangle} dy,$$
где $V^*$ --- сопряжённый конус. Её логарифм $\ln\varphi_V$ --- выпуклая функция, и каноническая риманова метрика Винберга--Кошуля $d^2\ln\varphi_V$ на конусе $V$ является метрикой Фишера--Рао при  интерпретации конуса $V^*$ как  многообразия вероятностных распределений.

Информационная геометрия является бурно развивающейся  областью математики и  применяется в самых различных областях прикладной математики, включая  нейронауки. Последнее связано с тем, что  мозг является машиной по переработке информации,  и  сформулированный  в рамках информационной геометрии  вариационный принцип К.~Фристона  
минимизации свободной энергии  предлагает объяснение  цели работы  мозга.

Информационная геометрия тесно связана с  геометрической термодинамикой, в частности, с термодинамикой Сурьо на группах Ли и  её  обобщениями. Согласно Ф.~Барбареско~\cite{Bar}, преобразование Лежандра от $\ln\varphi_V$  интерпретируется  в этой теории как энтропия.

Вклад  Э.~Б.~Винберга в развитие  информационной геометрии по достоинству был  оценён экспертами в этой области науки.  
Говоря о других применениях, стоит отметить, что в выпуклом программировании характеристическая функция $\varphi_V$ используется как барьерная  функция, позволяющая  для нахождения  минимума выпуклой функции в  выпуклом конусе применить классический метод Ньютона 
Теория $T$-алгебр широко используется в современной математике и теоретической физике, например, для  унификации  и обобщения   магического квадрата Фрейденталя (см.\ ниже) и  магической звезды Мукаи  
и при изучении чёрных дыр   в  пятимерной  $N=2$   супергравитации. 

\subsubsection*{Группы отражений в пространствах Лобачевского}

В середине 1960-х~гг.\ параллельно с теорией однородных кэлеровых многообразий Эрнест Борисович начал заниматься дискретными группами, порождёнными отражениями в гиперплоскостях в пространстве Лобачевского $\Lambda^n$ (гиперболическими калейдоскопами, как он любил их называть). В 1967~г.\ в <<Математическом сборнике>> вышла его большая статья \cite{Vinb67}. Эта работа состоит из трёх частей. В первой, геометрической, части работы строится систематическая теория гиперболических групп отражений в терминах матриц Грама фундаментальных многогранников, обобщающая классическую евклидову теорию Кокстера. Во второй, алгебраической, части доказывается критерий арифметичности групп отражений. В третьей части Винберг строит примеры арифметических и неарифметических групп отражений в~$\Lambda^n$. В частности, здесь естественным образом получены известные ранее примеры В.~С.~Макарова. 

Общая теория арифметических подгрупп в полупростых группах Ли была построена А.~Борелем и Хариш-Чандрой в~\cite{BorHC62}. Рассмотрим связную некомпактную полупростую группу Ли $\mathcal{G}$ с тривиальным центром. Пусть $G$ --- алгебраическая группа над вполне вещественным полем~$K$, у которой группа вещественных точек $G_{\RR}$ для некоторого вложения $K\hookrightarrow\RR$ локально изоморфна~$\mathcal{G}$, а для любого другого вложения эта группа 
компактна. Возьмём точное матричное представление $\rho$ группы $G$ и рассмотрим подгруппу
$$ G_I = \{ g\in G_{\RR} : \rho(g)_{ij} \in I \},$$
где $I$ --- кольцо целых поля~$K$. Подгруппы в~$\mathcal{G}$, соизмеримые с образами подгрупп вида~$G_I$, называются \emph{арифметическими}.

Данная конструкция является единственным известным общим методом построения решёток в~$\mathcal{G}$, т.е.\ дискретных подгрупп $\Gamma \subset \mathcal{G}$ с факторпространством $\mathcal{G}/\Gamma$ конечного объёма. В некоторых группах Ли вещественного ранга~1, например, в группах движений пространств Лобачевского существуют также неарифметические решётки. В работе \cite{Vinb67} Винберг получил критерий арифметичности групп отражений в терминах элементов матрицы Грама внешних нормалей к гиперграням фундаментального многогранника. Здесь же был определён класс квазиарифметических подгрупп, которые пользуются значительным интересом в последние годы. В заключительной части работы критерий арифметичности применяется для анализа конкретных примеров. На протяжении более пятидесяти лет фундаментальная и яркая статья Винберга является одной из основных ссылок в данной области.

В 1972 году в <<Математическом сборнике>> вышла работа Винберга~\cite{Vinb72}. В ней описывается получивший впоследствии широкую известность \emph{алгоритм Винберга} для построения арифметических групп, порождённых отражениями.  Примеры, полученные с помощью алгоритма в \cite{Vinb72} и в последующей совместной работе с И.~М.~Каплинской~\cite{VKa78}, имеют фундаментальное значение. Они нашли применение в работах Дж.~Конвея и Р.~Борчердса по гипотезе  <<moonshine>> (см., например, книгу Конвея и Слоана~\cite{Conway-Sloane}, одна из глав которой посвящена алгоритму Винберга и его группам).

В начале 1980-х годов Эрнест Борисович доказал глубокий и неожиданный результат об отсутствии решёток, порождённых отражениями, в пространствах $\Lambda^n$ большой размерности:

\begin{theorem}[\cite{Vinb81, Vinb84}]
В пространствах Лобачевского размерности $n \ge 30$ не существует дискретных групп отражений с ограниченным фундаментальным многогранником и арифметических дискретных групп отражений.
\end{theorem}

Для сравнения, в пространствах неотрицательной кривизны дискретные группы с ограниченным фундаментальным многогранником, порождённые отражениями, легко построить в любой размерности. Впоследствии, развивая некоторые идеи Винберга, М.~Н.~Прохоров доказал, что для $n \ge 996$ в $\Lambda^n$ не существует также и дискретных групп отражений с неограниченным фундаментальным многогранником конечного объёма. Группы отражений в пространствах Лобачевского и теоремы конечности были представлены Винбергом в его докладе на Международном конгрессе математиков в Варшаве~\cite{Vinb84_icm}.

Что происходит в тех размерностях, в которых существуют решётки и арифметические группы, порождённые отражениями? В размерности 2 полная классификация рефлективных решёток принадлежит Ф.~Клейну и А.~Пуанкаре. В размерности 3 классификация была получена Е.~М.~Андреевым, когда он учился в аспирантуре под руководством Э.~Б.~Винберга (1970~г.). Работы Андреева впоследствии легли в основу теории геометризации трёхмерных многообразий В.~П.~Тёрстона.  Систематическое изучение вопросов классификации арифметических групп отражений было начато в 1980-х~гг.\ в работах В.~В.~Никулина и Э.~Б.~Винберга, и затем продолжено в работах И.~Агола, М.~В.~Белолипецкого, К.~Маклахлана, А.~В.~Рида и других.

Э.~Б.~Винберг продолжал заниматься группами отражений в пространствах Лобачевского на протяжении всей математической карьеры. Его статья \cite{Vinb85} содержит изложение теории и обзор известных к тому времени результатов. Среди последующих результатов Винберга по этой теме можно упомянуть работу по подгруппам отражений в группах Бьянки~\cite{Vinb90}, совместную с Л.~Д.~Потягайло статью о прямоугольных группах отражений в $\Lambda^n$ \cite{PV05}, и также совсем недавнюю работу \cite{Vinb15} о неарифметических группах, порождённых отражениями, строящихся с помощью вариации метода Громова и Пятецкого-Шапиро. Как и другие работы Винберга, каждая из этих статей открывает новые горизонты для дальнейших исследований. 


\subsubsection*{Алгебраическая теория инвариантов}

Задачи теории инвариантов в широком смысле привлекали Эрнеста Борисовича со студенческих лет. Уже его дипломная работа, опубликованная в~\cite{V60} и посвящённая инвариантным линейным связностям в однородных пространствах, содержит классификацию локально-транзитивных неприводимых представлений простых групп Ли.

В совместной работе Э.~Б.~Винберга с учениками~\cite{AVE67} исследуется вопрос о три\-ви\-аль\-нос\-ти стационарной подалгебры общего положения (т.\,е.\ алгебры Ли стабилизатора вектора общего положения) для линейного представления комплексной полупростой группы Ли~$G$. Полученное достаточное условие
использует версию {\it индекса Дынкина\/} для представлений простой алгебраической группы.
Пусть $G\subset\GL(V)$ --- редуктивная комплексная линейная группа и $\g$ --- касательная алгебра Ли группы~$G$. Тогда билинейная форма на~$\g$, заданная формулой
$(x,y)\mapsto \tr_V (xy)$, невырождена. Если группа $G$ проста (как группа Ли), то отношение
$l_G=\tr_V (x^2)/ \tr_\g (\ad x)^2$ --- рациональное число, не зависящее от выбора
$x\in\g$; оно называется \emph{индексом} линейной группы~$G$. В \cite{AVE67} обнаружена неожиданная и удивительная связь индекса со стационарной подалгеброй общего положения (с.\,п.\,о.\,п.):

\begin{theorem}
Если $G$ полупроста, и для всех простых нормальных делителей $G_i$ группы $G$
верно, что $l_{G_i}>1$, то  
с.\,п.\,о.\,п.\ тривиальна.
\end{theorem}


Развивая методы этой работы, А.~Г.\ Элашвили получил вскоре классификацию всех представлений простых групп и неприводимых представлений полупростых групп с нетривиальной с.\,п.\,о.\,п.

Свойства нетривиальности стационарных подалгебр и стационарных подгрупп (стабилизаторов) векторов общего положения коррелируют с другими <<хорошими>> теоретико-инвариантными свойствами линейных представлений: алгебраической независимостью базисных инвариантных многочленов (порождающих алгебру инвариантов представления) и конечностью числа орбит в любом многообразии уровня инвариантных многочленов (где инварианты принимают заданные значения). Таких <<хороших>> представлений сравнительно немного, но они наиболее интересны с точки зрения теории инвариантов и её многочисленных приложений. В работе \cite{KPV76} Э.~Б.~Винберг, В.~Г.~Кац и В.~Л.~Попов получили список всех неприводимых представлений простых комплексных групп Ли со свободной 
алгеброй инвариантов, который совпал с ранее найденными списками линейных представлений в указанном классе, имеющих нетривиальный стабилизатор общего положения или конечное число орбит в любом многообразии уровня инвариантов.

Вопросы о причинах <<хорошего>> теоретико-инвариантного поведения линейных представлений и естественном источнике таких представлений давно интересовали исследователей. В своей фундаментальной работе 1963 года \cite{Kos} Б.~Костант построил теорию инвариантов присоединённого представления полупростой комплексной группы Ли $G$ в своей касательной алгебре Ли~$\g$. В частности, обнаружилось, что присоединённое представление обладает указанными выше <<хорошими>> свойствами. Рассказывая с восхищением о работе Костанта на своём семинаре, Е.~Б.~Дынкин сожалел лишь о том, что эти замечательные результаты не были получены его участниками, заметив, что <<Эрик (Винберг) вполне мог бы это сделать>>. Через десяток лет Эрнест Борисович оправдал интуицию своего учителя, создав новое направление в теории инвариантов, далеко развивающее и обобщающее результаты Костанта --- теорию тэта-групп.

Пусть полупростая алгебра Ли $\g$ снабжена периодической градуировкой по модулю~$m$:
$$
\g=\g_0\oplus\g_1\oplus\dots\oplus\g_{m-1}.
$$
Эта градуировка задаётся периодическим автоморфизмом $\theta$ алгебры~$\g$. Связная подгруппа Ли $G_0\subset G$, соответствующая подалгебре Ли $\g_0\subset\g$, действует в пространстве $\g_1$ путём ограничения присоединённого представления группы~$G$. Образ группы $G_0$ при этом представлении называется \emph{$\theta$-группой} градуированной алгебры~$\g$.

В 1973--1975~гг.\ Э.~Б.~Винберг доказал, что алгебра $G_0$-инвариантных многочленов на пространстве $\g_1$ свободна 
(это утверждение в виде гипотезы высказывал ещё в 1968~г.\ В.~Г.~Кац) и описал замкнутые $G_0$-орбиты в пространстве $\g_1$ \cite{D221,I40}.
В то же время Эрнестом Борисовичем был предложен метод носителей \cite{D225,S79} для классификации нильпотентных $G_0$-орбит в пространстве~$\g_1$. Вместе с описанием замкнутых (полупростых) орбит и градуированной версией разложения Жордана в полупростой алгебре Ли, этот результат позволяет в принципе описать все орбиты группы $G_0$ в пространстве~$\g_1$, что является конечной целью теории инвариантов для данного линейного представления. Испытательным полигоном для этого метода послужила полученная Э.~Б.~Винбергом
и А.~Г.~Элашвили классификация тривекторов 9-мерного пространства \cite{VE78}, что продолжило классические исследования В.~Райхеля, Я.~А.~Схоутена и Г.~Б.~Гуревича.

Теория тэта-групп работает также для непериодических (целочисленных) градуировок
$$
\g=\dots\oplus\g_{-2}\oplus\g_{-1}\oplus\g_0\oplus\g_1\oplus\g_2\oplus\cdots
$$
В этом случае группа $G_0$ не имеет непостоянных инвариантов в пространстве $\g_1$ и имеет лишь конечное число (нильпотентных) $G_0$-орбит. Можно рассмотреть естественную подгруппу $\widetilde{G}_0\subset G_0$ коразмерности~1 (или её образ в $\GL(\g_1)$ --- \emph{приведённую тэта-группу}), которая уже может иметь непостоянные инварианты, а её орбиты --- те же, что у~$G_0$, за исключением открытой плотной $G_0$-орбиты, которая может распадаться в объединение однопараметрического семейства замкнутых $\widetilde{G}_0$-орбит. 

Теория тэта-групп оказалась чрезвычайно полезной не только в теории инвариантов, теории представлений и теории алгебраических групп, но и в теории чисел, и в алгебраической геометрии, и теперь заслуженно называется теорией тэта-групп Винберга. Эта теория получила дальнейшее развитие над полями положительной характеристики~\cite{T12}.

Тэта-группы и приведённые тэта-группы почти исчерпывают список <<хороших>> неприводимых представлений простых групп Ли, о котором говорилось выше. Наиболее интересный из оставшихся случаев --- спинорное представление группы $\Spin_{13}$ --- был полностью исследован в работе Э.~Б.~Винберга (E.~Viniberghi) и В.~Г.~Каца (V.~Gatti) \cite{GaVi}. По совету Ж.-К.~Рота авторы взяли себе псевдонимы: Katze --- <<кошка>> по-немецки, а gatta --- по-итальянски; Viniberghi --- итальянизированное <<Винберг>>. В то время публикация статьи в соавторстве с Кацем, к тому времени эмигрировавшим в США, могла навредить Винбергу.

Говоря о работах Эрнеста Борисовича по геометрической теории инвариантов, нельзя не отметить и две
замечательные статьи, которые не так хорошо известны, как они того заслуживают, и
в которых ярко проявилось его умение найти новый подход к, казалось бы, уже хорошо
исследованным проблемам и минимальными средствами получить новые глубокие
результаты.

Используя метод носителей и критерий Гильберта--Мамфорда, Винберг получил в
работе~\cite{V80}  явное описание единственной замкнутой орбиты редуктивной линейной
группы $G\subset\GL(V)$,  содержащейся в замыкании произвольной $G$-орбиты
$G{\cdot}v\subset V$. Доказательство использует максимальную компактную подгруппу группы
$G$ и содержит немало полезных фактов о структуре замыкания $G$-орбит.
В качестве <<побочного продукта>> получается новый полезный критерий замкнутости $G$-орбиты.

Работа \cite{V00} содержит чисто алгебраический и, можно даже сказать, <<элементарный>>
подход к геометрическому понятию стабильности (т.\,е.\ замкнутости орбит общего положения) для действия редуктивной алгебраической группы $G$ на аффинном многообразии~$X$. Не используя почти ничего, кроме оператора Рейнольдса $f\mapsto f^\natural$ (проектора алгебры многочленов на подалгебру инвариантов) и понятия радикала идеала,
Эрнест Борисович строит идеал подмногообразия в~$X$, являющегося замыканием множества
всех замкнутых $G$-орбит в $X$ и получает критерий стабильности действия $G$ на $X$ как условие
невырожденности <<скалярного произведения>> $(f,g)\to (fg)^\natural$ на координатной алгебре многообразия~$X$. Из этого критерия мгновенно
получаются известные ранее результаты о стабильных действиях, а также и ряд новых (например, сохранение стабильности при ограничении действия на редуктивную подгруппу).

\subsubsection*{Алгебраические группы преобразований}

Алгебраическая теория инвариантов является частью более общей науки --- теории алгебраических групп преобразований, изучающей действия алгебраических групп на алгебраических многообразиях и их влияние на геометрию многообразий и теорию представлений. Эрнест Борисович начал работать в этой области с 1970-х годов.

В совместной работе Э.~Б.~Винберга с Б.~Н.~Кимельфельдом \cite{VKi78} получен критерий локальной транзитивности действия произвольной алгебраической подгруппы $H$ полупростой комплексной группы Ли $G$
на флаговом многообразии $G/P$. 
В частном случае, когда $P=B$ --- это борелевская подгруппа,
результат Винберга--Кимельфельда выглядит следующим образом:
\begin{theorem}
Действие $H$ на $G/B$ имеет открытую плотную орбиту тогда и только тогда, когда представление
группы $G$ в пространстве сечений любого однородного линейного расслоения $\LL\to G/H$ имеет простой спектр.
\end{theorem}
Вместо действия $H$ на $G/B$ можно рассматривать действие $B$ на $G/H$, и тогда предыдущая теорема даёт теоретико-представленческий критерий наличия открытой $B$-орбиты на однородном пространстве $G/H$. Однородные пространства или, более общо, алгебраические многообразия, снабжённые действием редуктивной группы~$G$, для которого борелевская подгруппа $B\subset G$ имеет открытую орбиту, получили название \emph{сферических}. Теорема Винберга--Кимельфельда --- один из краеугольных камней теории сферических многообразий, которая интенсивно развивается по сей день. Вот её версия, доказанная в той же работе:
\begin{theorem}
Если неприводимое $G$-многообразие $X$ сферично, то спектр представления группы $G$ в алгебре регулярных функций $\CC[X]$ прост. Для аффинных многообразий верно и обратное утверждение.
\end{theorem}

Говоря о вкладе Э.~Б.~Винберга в теорию сферических многообразий, нельзя не упомянуть его совместную работу с В.~Л.~Поповым \cite{VP72}. В ней для изучения одного класса сферических многообразий (которые тогда ещё так не назывались) впервые развита комбинаторно-геометрическая техника полиэдральных конусов и решёток, ставшая впоследствии одним из основных инструментов изучения сферических многообразий и их частного случая --- торических многообразий --- теория которых только зарождалась в те годы (начало 1970-х) в работах М.~Демазюра и др.

В работе \cite{V86} Эрнест Борисович ввёл понятие сложности действия
редуктивной алгебраической группы на неприводимом многообразии. Если $G,B,X$ --- такие же, как выше, то сложностью действия $G$ на $X$ называется коразмерность $B$-орбит общего положения в~$X$. 
Винберг получил замечательный результат:
\begin{theorem}[\cite{V86}]
Если $Y$ --- произвольное неприводимое $B$-инвариантное подмногообразие в~$X$, то
коразмерность $B$-орбит общего положения в $Y$ не превосходит сложности действия $G$ на~$X$.
\end{theorem}
Таким образом, наиболее <<богатое>> непрерывное семейство борелевских орбит состоит из орбит общего положения. Следует отметить, что, хотя утверждение касается действия группы~$B$, наличие действия объемлющей редуктивной группы $G$ существенно: без этого условия утверждение перестаёт быть верным.
Для действий сложности 0 (т.е.\ сферических) получается такое важное
\begin{corollary}
Если группа $B$ действует на $G$-многообразии $X$ с открытой плотной орбитой, то она имеет лишь конечное число орбит в~$X$.
\end{corollary}
Последнее утверждение независимо доказал М.~Брион. Оно появилось ещё в тексте работы \cite{VKi78} из-за опечатки, и уточняющий вопрос от читателя привёл Эрнеста Борисовича к результатам работы \cite{V86} (см.\ выше рассказ об этом).


Сложность оказалась тем инвариантом действия редуктивной группы, который отделяет <<ручные>> (поддающиеся исчерпывающему описанию) действия от <<диких>> (полная теория которых вряд ли возможна). О случае нулевой сложности (сферические многообразия) было сказано выше. К настоящему времени получена полная классификация сферических однородных пространств и их эквивариантных открытых вложений. Что касается действий сложности~$1$, то есть две возможности: либо $G$-орбиты общего положения образуют однопараметрическое семейство сферических однородных пространств, либо имеется открытая $G$-орбита сложности~$1$. Во втором случае возникает задача классификации всех однородных пространств группы $G$ сложности $1$ и их эквивариантных вложений. По этому поводу Э.~Б.~Винберг замечает в~\cite{V86}, что <<насколько мне известно, классификацией однородных пространств сложности $1$ никто не занимался, но кажется, что их должно быть немного>>.

В дальнейшем это предположение Эрнеста Борисовича полностью оправдалось. Классификация аффинных однородных пространств сложности $1$ для простых групп $G$ была получена Д.~И.~Панюшевым (1992~г.) и при этом оказалось, что половину достаточно короткого списка составляют подгруппы~$H$, получающиеся из сферических отбрасыванием одномерного центра. Полную классификацию аффинных однородных пространств сложности $1$ получили И.~В.~Аржанцев и О.~В.~Чувашова (2004~г.). Кроме того, замечено, что во многих классах 
$G$-многообразий (нильпотентные орбиты, аффинные двойные конусы и др.) многообразий сложности 1 значительно меньше, чем 
сферических. Общая теория вложений однородных пространств сложности $1$ была построена Д.~А.~Тимашевым (1997~г.).

В середине 1990-х~гг.\ Э.~Б.~Винберг обращается к теории алгебраических полугрупп. Алгебраическая полугруппа --- это неприводимое алгебраическое многообразие $X$ с ассоциативным умножением $X\times X\to X$ и единицей. Группа $G$ обратимых элементов алгебраической полугруппы $X$ --- это плотное открытое подмножество со структурой алгебраической группы. Теория алгебраических полугрупп разрабатывалась с начала 1980-х~гг.\ (М.~Путча, Л.~Реннер, и др.), но в работе \cite{V95} Винберг осуществил настоящий прорыв, получив полное описание \emph{редуктивных} алгебраических полугрупп (т.е.\ c редуктивной группой~$G$) и детально исследовав их алгебраические и геометрические свойства. Редуктивную алгебраическую полугруппу $X$ можно рассматривать как сферическое многообразие по отношению к действию группы $G\times G$ умножениями слева и справа, и идеи теории сферических многообразий активно используются в~\cite{V95}. Дальнейшее развитие теория редуктивных алгебраических полугрупп получила в работах А.~Риттаторе и др.

\subsubsection*{Пуассонова алгебра и симплектическая геометрия}

Статья \cite{V90} представляет собой первое обращение Эрнеста Борисовича к изучению пуассоновой структуры, связанной с коприсоединённым представлением группы Ли, и эквивариантной симплектической геометрии. В работе делается попытка квантования коммутативных (по Пуассону) подалгебр в симметрической алгебре $\Sym(\g)$ редуктивной алгебры Ли~$\g$, получаемых методом сдвига аргумента, восходящим к А.~С.~Мищенко и А.~Т.~Фоменко. Это сделано для элементов степени $2$ с помощью явной конструкции, которая потом использовалось другими авторами и для более общих целей. Подалгебры, получаемые методом сдвига аргумента, названы в работе
подалгебрами Мищенко--Фоменко (далее \emph{МФ-подалгебры}), и это название прижилось.
В работе также содержится несколько общих утверждений о связи коммутативных подалгебр в универсальной обёртывающей алгебре
$\Uni(\g)$ и в $\Sym(\g)$ и ещё скрытая жемчужина ---
явное описание выпуклой оболочки орбиты группы Вейля в пространстве весов, 
то есть описание граней соответствующего многогранника. Кроме того, в \cite{V90} показано, что в случае $\g=\sgl_n$, подалгебра
Гельфанда--Цетлина является пределом МФ-подалгебр.
В дальнейшем эта тематика заняла значительное место в деятельности как самого Эрнеста Борисовича, так и его учеников (А.~А.~Тарасов, В.~В.~Шувалов, О.~С.~Якимова, Л.~Г.~Рыбников, А.~А.~Зорин).

Статья~\cite{V90} реанимировала интерес к МФ-подалгебрам, привлекла внимание к задаче их
квантования и показала, что пределы МФ-подал\-гебр, пожалуй, даже интереснее их самих.
В 1996~г.\ М.~Л.~Назаров и Г.~И.~Ольшанский проквантовали МФ-подалгебры для
классических алгебр Ли при помощи янгианов, затем простая конструкция для
$\g=\sgl_n$ была предложена А.~А.~Тарасовым (2000~г.), а общее решение было дано в 2006~г.\
Л.~Г.~Рыбниковым, который использовал центр Фейгина--Френкеля, ассоциированный с
аффинной алгеброй Каца--Муди $\widehat{\g}$ и модель Годена. Дальнейшее изучение
связей между МФ-подалгебрами (равно как и их пределами) и подалгебрами Годена оказалось
необыкновенно плодотворно для обеих теорий.

В известной статье 1956~г.\ о формуле следа, А.~Сельберг ввёл \emph{слабо симметрические} римановы многообразия, которые впоследствии были основательно забыты.
Однако в 1999~г.\ появилась совместная работа Э.~Б.~Винберга с Д.~Н.~Ахиезером \cite{AV99},
в которой 
определяются слабо симметрические аффинные однородные многообразия редуктивной алгебраической группы $G$ над любым алгебраически замкнутым полем нулевой характеристики.
В случае поля комплексных чисел
получаются комплексификации компактных слабо симметрических римановых 
многообразий Сельберга. 

Основной результат работы состоит в том, что аффинное однородное пространство $X=G/H$ слабо симметрично тогда и только тогда, когда оно сферично. 
Одно из приложений состоит в том, что однородное риманово многообразие $Y$
слабо симметрично относительно связной \emph{вещественной} редуктивной группы $G$
тогда и только тогда, когда алгебра $G$-инвариантных дифференциальных
операторов на $Y$ коммутативна, т.е.\ когда $Y$ является \emph{коммутативным пространством}.

Две 
основополагающие статьи Э.~Б.~Винберга \cite{V01,V03} посвящены коммутативным однородным пространствам вещественных групп Ли (не обязательно редуктивных). Если
$Y=G/K$ --- риманово однородное пространство вещественной группы~$G$, то
коммутативность алгебры $G$-ин\-ва\-ри\-ант\-ных дифференциальных операторов на~$Y$ --- это одно из эквивалентных определяющих свойств \emph{пары Гельфанда} $(G,K)$ в классическом её понимании. Эрнест Борисович доказал, что в случае коммутативного однородного пространства группа $G$ 
разлагается в полупрямое произведение подгрупп Ли $G=L\ltimes N$, где $L$ редуктивна, а $N$ унипотентна. Более того, $K\subset L$, алгебра Ли $\n$ группы $N$ не более чем двуступенно нильпотентна, и всякий $K$-инвариантный многочлен на пространстве $\n$ $L$-инвариантен. Эти структурные результаты позволили Винбергу классифицировать неприводимые нильпотентные коммутативные пространства (т.е.\ такие, для которых $L=K$ и действие $K$ на
$\n/[\n,\n]$ неприводимо). В дальнейшем они привели к полной классификации пар Гельфанда в диссертации О.~С.~Якимовой (2005~г.), выполненной под руководством Эрнеста Борисовича.

Если однородное пространство $Y=G/K$ коммутативно, то действие $G$ на $T^*Y$
{\it коизотропно} (в смысле стандартной симплектической структуры кокасательного
расслоения). На момент написания статьи \cite{V01} обратная импликация была хорошо известна в случае редуктивной группы~$G$. Э.~Б.~Винберг показал, что она верна и в случае $L=K$ (\emph{пространства гейзенбергова типа}). Позднее Л.~Г.~Рыбников доказал эквивалентность в общем случае (2003~г.).

Коизотропные действия, а также вопросы связанные с {корангом} симплектических действий
относятся к предмету {эквивариантной симплектической геометрии}. В качестве приложений
здесь получаются утверждения об интегрируемых гамильтоновых системах.
В случае комплексной редуктивной группы~$G$, действующей на алгебраическом многообразии~$X$, коранг естественного действия $G$ на кокасательном расслоении $T^*X$ равен удвоенной сложности действия $G$ на~$X$. В частности, действие $G$ на $T^*X$ коизотропно тогда и только тогда, когда действие $G$ на $X$ сферично. В работе  \cite{V01} дано новое доказательство этого замечательного факта,
основанное на применении \emph{орисфер} (т.е. орбит максимальных унипотентных подгрупп группы~$G$).

Детальное изложение метода орисфер содержится в работе \cite{V01b}.
Там доказано, что если $X$ --- квазиаффинное многообразие с действием редуктивной группы
$G$ и $\Hor X$ --- многообразие орисфер общего положения в~$X$, то имеется каноническое $G$-эквивариантное симплектическое накрытие Галуа $T^*(\Hor X)\to T^*X$. Для (квазиаффинного) сферического однородного пространства $X$ многообразие $\Hor X$ также является сферическим однородным пространством.
Опираясь на этот факт, Эрнест Борисович последние годы работал над <<теорией Галуа>> сферических однородных пространств, в которой аналогом группы Галуа выступает группа вышеупомянутого накрытия. Эта работа осталась незавершённой.

\subsubsection*{<<Трансцендентная теория инвариантов>> и автоморфные формы}

В 2009 году Э.~Б.~Винберг возвращается к однородным симметрическим областям в связи с задачей о свободных алгебрах автоморфных форм
на эрмитовых симметрических пространствах некомпактного типа.

Пусть $X$ --- такое пространство, $G$ --- группа аналитических автоморфизмов пространства~$X$, а $\Gamma\subset G$ --- дискретная подгруппа
с факторпространством конечного объёма. Алгебра $\Gamma$-автоморфных форм  (относительно канонического фактора автоморфности)
$A(\Gamma)$  является неотрицательно градуированной алгеброй конечного типа. Пару $(X,\Gamma)$ назовём \emph{парой Шевалле}, если алгебра
$A(\Gamma)$  свободна.

Задача состоит в том, чтобы перечислить все такие пары (с точностью до естественной эквивалентности). Постановка задачи обязана
классической теореме Шевалле о свободных алгебрах инвариантов конечных линейных групп, и её решение придаст этой теореме
законченный и совершенный вид.

Для того, чтобы пара $(X,\Gamma)$  была парой Шевалле, необходимо, чтобы группа  $\Gamma$ порождалась комплексными отражениями.
Это условие оставляет для пространства $X$ только три возможности: комплексный  шар, полидиск или область $D_n$ четвёртого типа по классификации Э.~Картана ($D_2$~--- это произведение двух верхних
полуплоскостей, $D_3$~--- 3-мерная  верхняя полуплоскость Зигеля рода~2). Резонансный пример пары Шевалле
$({D_3},\Sp_4(\mathbb{Z}))$ был найден Игузой почти 60 лет тому назад, и на протяжении многих лет спорадические примеры
пар  Шевалле  в  размерностях, не больших~4, встречались в статьях разных авторов.

Прорыв произошёл в 2010 году, когда появилась первая из двух  работ \cite{V10,V18} Э.~Б.~Винберга, в  которых  впервые были построены
примеры пар Шевалле вида ($D_n,\Gamma$) для всех размерностей $5 \le n \le 10$. А в 2017 году,  в совместной работе \cite{VSh} Э.~Б.~Винберга и
О.~В.~Шварцмана, было показано, что  пары
такого вида могут существовать только при  $n \le 10$. Приступить к классификации дискретных групп, действующих в областях~$D_n$, со
свободной алгеброй автоморфных форм
Эрнест Борисович уже не успел, и те, кто доведёт дело до конца, будут многим ему обязаны.

\subsubsection*{Неабелевы градуировки алгебр Ли}

Завершим обзор ещё одной темой, которая
интересовала Э.~Б.~Винберга в течение всей его математической жизни. В 1964~г. Эрнест Борисович рассказал на семинаре П.~К.~Рашевского
по векторному и тензорному анализу о новой конструкции особых простых алгебр Ли,
возникшей под влиянием идей Б.~А.~Розенфельда, см.~\cite{V66}.
Эта конструкция использует две композиционные алгебры $A_1$, $A_2$ (над полем
$\CC$ или~$\RR$) и пространство косоэрмитовых матриц третьего порядка с
элементами из алгебры ${A_1\otimes A_2}$.

Как известно, всякая композиционная алгебра
над полем характеристики $\ne 2$ имеет размерность 1, 2, 4 или~8.
Перебор всех возможных пар $(A_1, A_2)$ над $\CC$ даёт квадрат $4\times4$ (<<магический квадрат Фрейденталя>>), включающий в себя все особые простые алгебры Ли, кроме~$\Gg_2$. В отличие от сходной конструкции Ж.~Титса, полученной чуть ранее, конструкция Винберга симметрична относительно $A_1$ и $A_2$ и тем самым объясняет симметрию <<магического квадрата>>. Интереснейшее замечание, сделанное Эрнестом Борисовичем, говорит, что если вместо матриц 3-го порядка брать матрицы произвольного порядка и ограничиться ассоциативными композиционными алгебрами, то эта конструкция даёт все комплексные и вещественные классические алгебры Ли.

Сорок лет спустя, Винберг вернулся к этой тематике в работе~\cite{V05}, где он ввёл понятия коротких $\SL_2$-структур и триад в простой группе Ли $G$ и доказал с их помощью гипотезу Розенфельда о существовании эллиптических плоскостей над тензорным произведением двух композиционных алгебр. Он также доказал, что все симметрические пространства особых групп Ли являются квазиэллиптическими плоскостями над некоторыми неассоциативными алгебрами (не обязательно тензорными произведениями композиционных алгебр).

Систематизация этого подхода привела Эрнеста Борисовича к общему понятию неабелевой градуировки на алгебре Ли~$\g$. $S$-структурой на $\g$ называется действие редуктивной группы $S$ автоморфизмами алгебры~$\g$. $S$-структура задаёт разложение пространства $\g$ в прямую сумму изотипных компонент:
$$
\g\simeq\bigoplus_iV_i\otimes U_i,
$$
где $V_i$ --- различные неприводимые представления группы~$S$, а $U_i=\Hom_S(V_i,\g)$ --- пространства кратностей. Такие изотипные разложения Эрнест Борисович называл неабелевыми градуировками (обычные градуировки с помощью циклических групп получаются, если в качестве $S$ взять $\CC^{\times}$ или циклическую группу конечного порядка). Скобка Ли на $\g$ задаёт на пространствах $U_i$ разные интересные алгебраические структуры (йордановы алгебры, алгебры Кантора--Аллисона и др.)

В работе \cite{V17} Э.~Б.~Винберг исследовал короткие $\SL_3$-структуры на простых алгебрах Ли (в которых $\dim V_i\le3$), а в 2017 году он прочел курс лекций <<Неабелевы градуировки простых алгебр Ли>> на зимней школе по теории инвариантов в МГУ. Эрнест Борисович положил начало этому новому направлению, и тут остаётся много работы для его последователей.

\subsubsection*{Полный список публикаций Э.~Б.~Винберга на Math-Net.ru:}

\verb|http://www.mathnet.ru/php/person.phtml?personid=8541|\\
\strut\hfill\verb|&wshow=&showmode=simple&option_lang=rus#pubs|

\renewcommand\refname{Избранные публикации Э.~Б.~Винберга}

\renewcommand\refname{Другие цитируемые работы}

\end{document}